\documentclass[11pt]{article}
\usepackage{amssymb,latexsym}

\usepackage{amsmath,amscd}
\usepackage{theorem}
\usepackage[all]{xy}
\title{\bf The integration problem for complex Lie algebroids}
\author{Alan Weinstein\thanks{research
    partially supported 
by NSF grant DMS-0204100
\newline \mbox{~~~~}MSC2000 Subject Classification Numbers: 
\newline \mbox{~~~~}Keywords: }
 \\      Department of Mathematics\\
        University of California\\ Berkeley
        CA, 94720-3840 USA\\{\tt alanw@math.berkeley.edu}}

\newcommand{\frakg}      {\mathfrak{g}}
\newcommand{\frakk}      {\mathfrak{k}}

\newcommand{\del}	{\partial}

\newcommand{\cale}       {\mathcal{E}}

\newcommand{\calj}       {\mathcal{J}}

\newcommand{\calt}       {\mathcal{T}}

\newcommand{\vbar}{\overline{v}}

\newcommand {\pf}{\noindent{\bf Proof.}\ }
\newcommand{\qed}{\begin{flushright} $\Box$\ \ \ \ \ \end{flushright}}

\newcommand{\reals}      {{\mathbb R}}
\newcommand{\integers}   {{\mathbb Z}}
\newcommand{\complex}    {{\mathbb C}}

\newtheorem{thm}{Theorem}[section]

\newtheorem{lemma}[thm]{Lemma}

\newtheorem{dfn}[thm]{Definition}

\newtheorem{rmk}[thm]{Remark}
\newtheorem{ex}[thm]{Example}

\begin{document}

\maketitle

\begin{abstract}
A complex Lie algebroid is a complex vector bundle over a smooth
(real) manifold $M$ with a bracket on sections and an anchor to the
complexified tangent bundle of $M$ which satisfy the usual Lie
algebroid axioms.  A proposal is made here to integrate analytic complex Lie
algebroids by using analytic continuation to a complexification of $M$
and integration to a holomorphic groupoid.  A collection of diverse
examples reveal that the holomorphic stacks presented by these
groupoids tend to coincide with known objects associated to structures
in complex geometry.  This suggests that the object integrating a
complex Lie algebroid should be a holomorphic stack.

\end{abstract}

\section{Introduction}

It is a pleasure to dedicate this paper to Professor Hideki Omori.  His work over many years, introducing  ILH manifolds \cite{om-group}, Weyl manifolds
\cite{om-ma-yo:weyl}, and blurred Lie groups \cite{om-ma-mi-yo:strange} has broadened the notion of what constitutes
a ``space."   The problem of  ``integrating'' complex vector fields
on real manifolds seems to lead to yet another kind of space, which is investigated in this paper.  

Recall that a {\bf Lie algebroid} over a smooth manifold $M$ is a
real vector
bundle $E$ over $M$ with a Lie algebra
structure (over $\reals$) on its sections and a bundle map
$\rho$ (called the {\bf anchor}) from $E$ to the tangent
bundle $TM$, satisfying
the Leibniz rule
$$[a,fb]=f[a,b]+(\rho(a)f)b$$ for sections $a$ and $b$ and smooth
functions
 $f:M\to \reals$.  Sections of a Lie algebroid may be
thought of as ``virtual" vector fields, which are mapped to ordinary
vector fields by the anchor.
 
 There is an analogous definition for
complex manifolds, in which $E$
 is a holomorphic vector bundle over
$M$, and the Lie algebra structure
 is defined on the sheaf of local
sections.  Such objects are called
 complex Lie algebroids by Chemla
\cite{ch:duality}, but they will be called in this paper {\bf
holomorphic Lie algebroids} to distinguish them from the ``hybrid"
objects defined in \cite{ca-we:geometric} as follows.

\begin{dfn}
A {\bf complex Lie algebroid} ({\bf CLA}) over a smooth (real) manifold $M$
is a complex
 vector bundle $E$ over $M$ with a Lie algebra structure
(over $\complex$)
 on its space $\cale$ of sections and a bundle map
$\rho$ (called the {\bf
 anchor}) from $E$ to the complexified
tangent bundle $T_\complex M$, satisfying
 the Leibniz rule
$$[a,fb]=f[a,b]+(\rho(a)f)b$$ for sections $a$ and $b$ in $\cale$ and 
smooth functions $f:M\to \complex$
\end{dfn}

The unmodified term ``Lie algebroid''  will always mean
``real Lie algebroid."

Every Lie algebroid may be realized as the bundle whose
sections are the left invariant vector fields on a {\em local} Lie groupoid
$\Gamma$.  The {\bf integration problem} of determining when
$\Gamma$ can be taken to be a global groupoid was
completely solved in \cite{cr-fe:lie}, but, for a complex Lie algebroid,
it is not even clear what the corresponding local object should be.
 The main purpose of the present paper is to
propose a candidate for this object.  

Any  CLA $E$ whose anchor is injective  may be identified with
the involutive subbundle $\rho(E) \subseteq T_\complex M$.
Such subbundles have been studied extensively under the name of
``involutive structures'' or ``formally integrable structures,'' for
instance by Treves \cite{tr:hypo}.   An important issue in these studies
has been to establish the existence (or nonexistence) of ``enough integrals,"  i.e. smooth
functions which are annihilated by all the sections of $E$.  In the
general $C^\infty$ case, the question is very subtle and leads to
deep problems and results in linear PDE theory.  When $E$ is 
analytic,\footnote{In this paper, ``analytic'' will always mean ``real
analytic'', and ``holomorphic'' will be used for ``complex analytic."} though,
one can sometimes proceed in a fairly straightforward way by
complexifying $M$ and extending $E$ by holomorphic continuation to an
involutive holomorphic tangent subbundle of the complexification,
where it defines a holomorphic foliation.  The leaf space of this
foliation is then a complex manifold whose holomorphic functions
restrict to $M$ to give integrals of $E$.  (A succinct example of this
may be found at the end of \cite{ro:review}; see
Section \ref{subsubsec-cr} below.)  

The leaf space described above may be thought of as
the ``integration" of the involutive subbundle $E$; 
this suggests a similar approach to analytic CLAs whose anchors may
not be injective.  Any analytic CLA $E$ over $M$ may be 
holomorphically continued to a holomorphic Lie algebroid $E'$ over a
complexification $M_\complex$; $E'$ may then be integrated to a (possibly
local) holomorphic groupoid $G$.   Since $G$ will generally have
nontrivial isotropy, one must take this into account by considering
not just the orbit space of $G$, but the ``holomorphic stack" 
associated to $G$.

Some intuition behind the complexification approach to integration
comes from the following picture in the real case.  If $G$ is a Lie
group, there is a long tradition of thinking of its Lie algebra
elements as tiny arrows pointing from the identity of $G$ to
``infinitesimally nearby" elements.  If $G$ is now a Lie groupoid over
a manifold $M$, $M$ may be identified with the identity elements of $G$,
and an element $a$ of the Lie algebroid $E$ of $G$ may be thought of as
an arrow from the base $x\in M$ of $a$ to a groupoid element with
its source at $x$ and its target at an infinitesimally nearby 
$y \in G$.  The tangent vector $\rho(a)$ is then viewed as a tiny
arrow in $M$ pointing from $x$ to $y$.

Now suppose that $E$ is a complex Lie algebroid over $M$.  Then
$\rho(a)$ is a complex tangent vector.  To visualize it, one may still
think of the tail of the tiny arrow as being at $a$, but the imaginary
part of the vector will force the head to lie somewhere ``out there"
in a complex manifold $M_\complex$ containing $M$ as a totally real
submanifold.  To invert (and compose) such groupoid
elements requires that their sources as well as targets be allowed to lie in this
complexification $M_\complex$.  Thus, the integration should be a
groupoid over the complexification.

What exactly is this complexification?  Haefliger
\cite{ha:structures}, Shutrick \cite{sh:complex}, and Whitney and Bruhat 
 \cite{wh-br:quelques} all showed that every analytic
manifold $M$ may be embedded as an analytic, totally real submanifold of
a complex
manifold $M_\complex$.  Any two such complexifications are canonically
isomorphic near $M$.    Consequently, the identity
map extends uniquely near $M$ to an antiholomorphic involution of $M_\complex$
(``complex conjugation") having $M$ as its fixed point set.  Finally,
Grauert \cite{gr:levi} showed that
the complexification may be taken  to have a
pseudoconvex boundary and therefore be a Stein manifold.  $M_\complex$
is then called a {\bf Grauert tube}.

Of course, constructing the complexification requires that the Lie
algebroid have a real analytic structure.  For the underlying smooth manifold $M$, such
a structure exists and is unique up to isomorphism
\cite{wh:differentiable}, though the isomorphism between two such
structures is far from canonical.  Extending the analyticity to $E$ is
an issue which must be deal with in each example.  

In fact, examples are at the heart of this paper.  Except for some
brief final remarks, the many observations and questions about CLAs
which arise naturally by extension from the real theory and from complex
geometry will be left for future work.   Concepts such as
cohomology, connections, modular classes, K\"ahler structure, and
quantization are discussed by Block
\cite{bl:duality} and Cannas and the author \cite{ca-we:geometric} and
in work in progress with Eric Leichtnam and
Xiang Tang \cite{le-ta-we:poisson}.

\noindent
{\bf Acknowledgements.}  This ideas in this paper have developed over several years, in part during visits to MSRI,  Institut Math\'ematique de Jussieu, and \'Ecole Polytechnique.  I would like to thank these institutions for their support and hospitality, and many people for their helpful comments, including Marco Gualtieri, Yvette Kosmann-Schwarzbach, Claude Lebrun, Eric Leichtnam, Pierre Schapira, Xiang Tang, and Fran\c{c}ois Treves.  Finally, I thank Asha Weinstein for editorial advice.  

\section{Complexifications of real Lie algebroids}

A complex Lie algebroid over a point is just a Lie algebra $\frakg$ over
$\complex$.  It seems natural to take as integration of $\frakg$ a
holomorphic Lie group $G$ with this Lie algebra.  In particular, if
$\frakg$ is the complexification of a real Lie algebra
$\frakg_\reals$, then $G$ is a complexification of a real Lie group
$G_\reals$.

Next, given any real Lie algebroid $E_\reals$, its complexification $E$
becomes a complex Lie algebroid when the bracket and anchor are extended
by complex (bi)linearity.  If $E_\reals$
is integrated to a (possibly local) Lie groupoid $G_\reals$, then a natural
candidate for $G$ would be a complexification of
$G_\reals$.  For this complexification to exist, 
$G_\reals$ must admit an analytic structure, and, when this structure does exist, it
is rarely unique (though it may be unique up to isomorphism).   

\subsection{Zero Lie algebroids}

Let $E_\reals$ be the zero Lie algebroid over $M$.  An analytic
structure on $E_\reals$ is just an analytic structure on $M$, which
exists but is unique only up to noncanonical isomorphism.  Now the
unique source-connected Lie groupoid integrating $E_\reals$ is the
manifold $M$ itself, which always admits a
complexification $M_\complex$.  This complex manifold is far from unique,
but its germ along $M$ is unique up to natural (holomorphic)
isomorphism, given the analytic structure on $M$.  One could say that
the choice of analytic structure on $M$ is part of the integration of
this zero complex Lie algebroid.

This example suggests that the object integrating $M$ should be the germ along $M$ of a complexification of $M$.  
Getting rid of all the choices, including that of the
analytic structure, requires that the
complexification $M_\complex$ be shrunk even further, 
to a formal neighborhood of $M$ in
$M_\complex$.  Both of these possibilities will be considered in many
of the examples which follow.

\begin{rmk}
{\em One could define the germ as an object for which the underlying
topological space is $M$, but with a structure sheaf given by germs
along $M$ of holomorphic functions on $M_\complex$.  But these are just
the analytic functions on $M$.  For the formal
neighborhood, the structure sheaf becomes simply the infinite jets of
smooth complex-valued functions.}
\end{rmk}

\subsection{Tangent bundles}
Let $E=T_\complex M$ be the full complexified tangent bundle.  Once
again, an analytic structure on $E_\reals=TM$ is tantamount to an
analytic structure on $M$, which leads to many  complexifications
$M_\complex$, as above.  A source-connected Lie groupoid
integrating $TM$ is the pair groupoid $M\times M$, while the
source-simply connected groupoid is the fundamental groupoid $\pi(M)$.
The pair groupoid $M_\complex \times M_\complex$ is then a
complexification of $M \times M$ and may be taken as an integration of
the complex Lie algebroid $T_\complex M$.  On the other hand, 
$\pi(M)$ could be complexified to $\pi(M_\complex)$; however, the result is 
sensitive, even after restriction to $M$, to the choice of
$M_\complex$.  If $M_\complex$ is taken to be a small neighborhood
of $M$, the restriction to $M$ is just $\pi(M)$.

\subsubsection{Interlude: The integration as a stack}
Some of the dependence on the choice of $M_\complex$ disappears when
two  groupoids are declared to be ``the same" when they are Morita equivalent. The groupoid is then seen as a presentation of a differential stack (see Behrend \cite{be:cohomology}
and Tseng and Zhu \cite{ts-zh:integrating}) or, more precisely, a holomorphic stack.  Since a transitive groupoid is equivalent to any of its isotropy groups, the stack represented by a pair groupoid $M\times M$ is just a point.  The only difference between this and 
$M_\complex\times M_\complex$ is that the latter represents a ``holomorphic point."  Depending on the choice of groupoid, this point as a stack might carry isotropy equal to the fundamental group of $M$ or even of one its complexifications.  

\subsection{Action groupoids}
Any (right) action of a Lie algebra $\frakk$ on $M$ induces an 
action, or transformation, groupoid structure on the trivial vector
bundle $E_\reals =M\times \frakk$.  The complexified bundle
$E=M\times\frakk_\complex$ becomes a complex Lie algebroid whose
anchor maps the constant sections of $E$ to a finite dimensional Lie algebra of complex vector fields on $M$.

When the $\frakk$ action comes from a (left) action of a Lie group $K$, $E_\reals$ integrates to the transformation groupoid $K\times M$; in fact, Dazord  \cite{da:groupoide} showed that $E_\reals $ is always integrable to a global groupoid $G$ which encodes the (possibly local) integration of the $\frakk$ action. 

Passing from $E_\reals$ to $E$ complicates issues significantly.
First, complexifying $G$ requires an analytic structure
on it, which amounts to an analytic structure on $M$ for which
the $\frakk$ action is analytic.  But this can fail to exist even when 
when $\frakk=\reals$, in other words, when the action is simply given by a vector
field.  For instance, if the vector field vanishes to infinite order
at a point $p$ of $M$, but not on a neighborhood of $p$, it can never
be made analytic, so complexification of the action groupoid $G$ and
hence integration of $E$ become impossible except on the formal level.

In addition, it is conceivable that some smooth action groupoids may be made analytic in essentially different ways, even though, according to Kutzschebauch 
\cite{ku:uniqueness}, this cannot happen for proper actions by groups with finitely many connected components.  Perhaps there is a smooth actions which admits several quite different complexifications.

When $M$ and the $\frakk$ action are analytic, the
vector fields generating the action extend to holomorphic vector fields on a
complexification $M_\complex$, leading to a holomorphic Lie algebroid
structure on $M_\complex \times k_\complex$.  This integrates to a
holomorphic Lie groupoid $G$, the ``local transformation groupoid" of
the complexified $K_\complex$ action.

Note the slightly different strategy here--the Lie
algebroid is first extended to the complexification and then
integrated, 
rather than
the other way around.  This  strategy will be used
extensively below.

\begin{ex}
{\em Let $\frakk = \reals$ act on $M=\reals$ via the vector field $x\,
\frac{\del}{\del x}$.  When $\frakk$ is considered as the Lie algebra of the
multiplicative group $\reals^+$, the resulting action groupoid
is $\reals^+ \times \reals$, with the first component acting on the
second by multiplication.  The orbits of this groupoid are the two
open half lines and the origin.

A natural complexification of  $\reals^+ \times \reals$ is the action groupoid 
$\complex^\times \times\complex$, whose orbits are the origin in
$\complex$ and its complement $\complex ^\times$.    When this
groupoid is restricted to the original manifold $\reals$, the two half
lines now belong to the same orbit, even if the
complexification $\complex$ is replaced by a small neighborhood of the
real axis.  (In this case, the complexified groupoid would no longer
be an action groupoid, but it would still have just the two orbits.)
As a stack, the complexified groupoid represents a space with two points, one of which is an ordinary holomorphic point.  The second point is in the closure of the first and has isotropy group  $\complex^\times$.

After restriction of the groupoid to the germ of $\complex$ around
$M$, or to the formal neighborhood, 
the notion of ``orbit" is harder to pin down, since the groupoid
does not directly define an equivalence relation.  

A somewhat different result is obtained if the algebroid is first
extended and then integrated.  The extended complex Lie algebroid is $\complex \times \complex$; for its natural integration,  the group is the simply connected cover $\complex$ of $\complex^\times$.  The action groupoid is now $\complex \times \complex$ with the action $w\cdot z=e^w z$, for which the orbits are the same as before, but  the isotropy group of nonzero $z$ (including real $z$) is now $2\pi i \integers$.  
}
\end{ex}

\begin{rmk}
{\em
A similar but slightly more complicated example is given by the vector
field on the phase plane $M=\reals^2$ which describes a classical
mechanical system near a local maximum of the potential function.  The
complexication of the action groupoid $\reals \times \reals^2$ includes groupoid elements connecting states on opposite sides of the potential maximum which cannot be connected by real classical trajectories.   These groupoid elements are not without physical interest, though, since they may be interpreted as representing quantum tunneling.
}  
\end{rmk}

\subsection{Foliations}
An analytic foliation $E_\reals\subset TM$ extends to
a holomorphic foliation of $M_\complex$, and, if $M_\complex$
is small enough,  the leaf stack of the
latter is just a straightforward complexification of the (analytic)
leaf stack of the former.  In particular, if the former is a manifold,
so is the latter.  

But there are many foliations which admit no compatible analytic
structure.  Take for example the Reeb  \cite{re:certaines} foliation
(or for that matter, according to Haefliger  \cite{ha:structures},
any foliation) on $S^3$.  The leaf space of the
Reeb foliation consists of two circles and a special point whose only
open neighborhood is the entire space.  The isotropy group of the holonomy
groupoid is trivial  for the leaves on the circles and $\integers ^2$
for the special leaf.  

To complexify the Lie algebroid by complexifying the
foliation groupoid, one might look instead at the equivalent groupoid
given by restriction to a cross section to the leaves.  This cross
section can be taken as a copy of $\reals$ on which $\integers ^2$
acts, fixing the origin, with one of the two generators acting by
1-sided contractions on the left half line and the other by
contractions on the right.  Complexifying the action of
the  generators gives maps on $\complex$ which have essential
singularities at the origin, and there seems to be no way to make a
holomorphic stack out of this data. 

\section{Involutive structures}
A complex Lie algebroid $E$ over $M$ with injective anchor may be
identified with the image of its anchor, which is an involutive
subbundle of $T_\complex M$.  Following Treves \cite{tr:hypo}, these
subbundles will be called here {\bf involutive structures}.  An analytic structure on $E$ is just an analytic structure on $M$ for which $E$ admits local bases of analytic complex vector fields.  

Let $E$ be an analytic subbundle of $T_\complex M$, then, and
$M_\complex$ a complexification of $M$.  Identifying
$T_\complex M$ with the restriction to $M$ of $TM_\complex$, one may extend
the local bases of analytic
sections of $E$  to local holomorphic sections of
$TM_\complex$.  For $M_\complex$ sufficiently small, local
bases again determine a holomorphic subbundle $E'$ of
$TM_\complex$.  Holomorphic continuation of identities implies that $E'$
is itself involutive;  by the holomorphic Frobenius theorem, 
it determines a holomorphic foliation of
$M_\complex$.  The holonomy groupoid of this foliation determines a
holomorphic stack which may be considered as the integration of the complex Lie
algebroid $E$.

The rest of this section is devoted to examples of involutive structures viewed as CLAs.  

\subsection{Complex structures}
\label{subsubsec-complex}
An almost complex structure on $M$ is an endomorphism 
$J:TM\to TM$ such that $-J^2$ is the identity.  $T_\complex M$ is the
direct sum of the $-i$ and $+i$ eigenspaces of the complexified
operator $J_\complex$. These  conjugate complex subbundles, 
 denoted by
$T^{0,1}_J M$ and $T^{1,0}_J M$ respectively, are involutive if and only if $J$ is integrable in the
sense that the Nijenhuis tensor $N_J$ vanishes.  The eigenspace
$T^{0,1}_J M$ is then a CLA which, like $J$ itself, is called a {\bf
complex structure}.  It is a standard fact that every subbundle
$E\subset T_\complex M$ such that $E\oplus \overline E=T_\complex M$
is $T^{0,1}_J M$ for some almost complex structure $J$.

Theorems of Eckmann-Fr\"olicher 
\cite{ec-fr:integrabilite} and Ehresmann \cite{eh:varietes} (analytic case)\footnote{The cited authors also attribute the result to de Rham.} and 
Newlander-Nirenberg \cite{ne-ni:complex} (smooth case) tell
us that any complex structure on $M$ is locally isomorphic to the
standard one on $\reals^{2n}=\complex^n$; i.e., it gives a reduction of
the atlas of smooth charts on $M$ to a subatlas with holomorphic
transition functions, making $M$ into a complex manifold.  Let us
pretend for a moment, though, that we do not know those theorems and look
directly at the integration of an analytic complex structure as a
holomorphic stack.  (The result of this exercise will turn out to be the original
1951 proof!) 

According to the discussion above, complexification gives a foliation
$E'$ of a suitably small $M_\complex$ whose leaves, by the condition
$E\oplus \overline E=T_\complex M$, have tangent spaces along $M$
which are complementary to the real subbundle $TM$.  As a result,
shrinking $M_\complex$ again can insure that each leaf is a
holomorphic ball intersecting $M$ exactly once, transversely, so that
the leaf space of this foliation may be identified with $M$.  This
leaf space being a complex manifold, $M$ itself inherits the structure
of a complex manifold.  Holomorphic local coordinates on $M$  result
from sliding open sets in $M$ along the foliation $E'$ to identify
them with open sets in holomorphic transversals, e.g. leaves of the holomorphic foliation $\overline{E'}$ which extends $\overline{E}$.

The holomorphic stack in this case may be identified with $M$ as a complex manifold, presented by the holonomy groupoid of the foliation $E'$.  An alternate presentation is the etale groupoid obtained by restricting the holonomy groupoid to the union of enough transversals to cover $M$ under projection along $E'$.  The latter groupoid is just the equivalence relation associated to a covering of $M$ by holomorphic charts.

When $E$ is given simply as a smooth complex structure, the only
recourse is to invoke the Newlander-Nirenberg theorem.  This has the
consequence that $M$ has an analytic structure in which $E$ is
analytic, so the previous situation is obtained.

\begin{rmk}
{\em The analytic structure on $M$ which makes a complex structure $E$ analytic is unique, since it must be the one attached to the holomorphic structure determined by $E$.  The situation is therefore different from that for the complex Lie algebroid
 $T_\complex M$ and the zero Lie algebroid, whose integration depends on an arbitrary choice of analytic structure compatible with the given smooth structure.
}
\end{rmk}

\subsection{CR structures}
\label{subsubsec-cr}
A step beyond the complex structures within the class of involutive systems are the general {\bf CR structures}.  These are subbundles $E$ of $T_\complex M$ such that $E$ and 
$\overline{E}$ intersect only in the zero section, but $E\oplus \overline{E}$ is not necessarily all of $T_\complex M$.\footnote{Some authors use the term ``CR structure" only when $E\oplus \overline{E}$ is of codimension one in $T_\complex M$.}

 Any ``generic" real submanifold $M$ in a
complex manifold $X$ 
inherits a CR structure, namely the intersection
$G_{M,X}=T_\complex M \cap T^{0,1}_J X$.  To be precise, the submanifold is called generic when $G_{M,X}$ has constant dimension; note that real hypersurfaces are always generic in this sense.
$G_{M,X}\oplus \overline{G_{M,X}}$
is the complexification of
the maximal complex subbundle $F_{M,X}$ of $TM$.  
A natural geometric problem, which has led to fundamental developments
in linear PDE theory, is whether a given CR manifold can be realized either
locally or globally as a submanifold in some complex manifold, and
in particular in $\complex^n$.  For analytic CR structures, the integration method of this paper solves this problem.  What follows below  essentially reproduces an argument of Andreotti and Fredricks
\cite{an-fr:embeddability}, or more precisely, that in the review by Rossi \cite{ro:review} of that paper.

Let $E'$ be the integrable holomorphic subbundle of  $TM_\complex$
which extends $E$.  The corresponding foliation will be 
called the {\bf CR foliation}. If $M$ has (real) dimension $2n+r$ and $E$ has complex dimension $n$, then $M_\complex$ has
complex dimension $2n+r$, and the leaves of the CR foliation  have complex dimension $n$; each of them meets $M$ in a
point, with no common tangent vectors (since $E$ contains no real
vectors).  It follows that $M_\complex$ can be chosen 
so that the leaves are simply connected; the
stack defined by the foliation groupoid is then 
simply a complex manifold $N$ of complex dimension $n+r$ containing $M$ as a real
hypersurface of real codimension $r$.  When $r=0$, $N=M$, and $M$ is a complex manifold; when $n=0$ (zero Lie algebroid), $N=M_\complex$.  
(Andreotti and Fredricks \cite{an-fr:embeddability} call  $N$  a complexification of $M$ for any $n$; thus, the complexification of a complex manifold is the manifold itself.)

\subsection{The Mizohata structure}

The next example shows that the natural map from $M$ to a stack which integrates a complex Lie algebroid $E\to M$ may not be injective.

As in  Example I.10.1 of Treves \cite{tr:hypo}, the {\bf Mizohata structure} over 
$M=\reals^2$ is defined to be the involutive system $E$  spanned by the complex vector field 
$$i\del/\del t -t\ \del/\del x.$$  It is a complex structure except
along the $x$-axis, where it is the complexification of the real
subspace spanned by $\del/\del t$.  The holomorphic continuation of $E$
over $\complex ^2$ is spanned by the same vector field in which
$(x,t)$ are taken as complex variables, and the leaves of the
corresponding foliation $E'$ are the levels of the invariant function
$\zeta =x-it^2/2$.  These levels, which can be described as graphs 
$x=it^2/2+\zeta$ with the parameter $t$ running through $\complex$, are
contractible, so the stack defined by the foliation groupoid is
isomorphic to $\complex$ with $\zeta$ as its complex coordinate.  The
natural map from $M$ to this stack folds $\reals^2$ along the $x$-axis, and the image is the (closed)
lower half plane.

The situation becomes more complicated rather than simpler if
the complexification is shrunk to a neighborhood of $\reals^2$ in
$\complex^2$, for instance that defined by the bounds $|\Im t| < \epsilon$ and
$|\Im z| < \epsilon$ on the imaginary parts.  
In this case, some of the level manifolds of $\zeta$
split into two components, so that
the corresponding part of the leaf space (the complement of a strip
near the origin in the lower half plane) bifurcates into two
branches.\footnote{There is no bifurcation in the upper half plane.}
The common closure of these branches is a family
of leaves depending on one (real) parameter, so we
can describe the integration of the Mizohata structure (or the
``complexification'', in the language used in CR geometry) as the
non-Hausdorff complex manifold which is the 
union of an open strip along the real axis in the complex
$\zeta$-plane with two copies of the rest of the lower half
plane.  The map from $M$ to this stack now separates points except
those in a strip along the $x$ axis, which is folded as before.

Integrals of the involutive structure on $M$ must
be even in $t$ near the $x$ axis; since they are holomorphic away from the
$x$ axis, they must be even everywhere.  In this case,
there are integrals of $E$ which are not the pullback of holomorphic
functions on the stack.  (See Example III.2.1 in Treves
\cite{tr:hypo}.)    

It is not clear what kind of geometric object is
obtained in the limit as the complexification shrinks
down to $M$, or for the formal complexification.

A test problem for any global theory of integration is to describe the integration of involutive structures on smooth surfaces which have singularities along a collection of simple closed curves and which are complex structures elsewhere.

\subsection{Eastwood-Graham and LeBrun-Mason structures}

In the next example, due to Eastwood and Graham \cite{ea-gr:involutive}, the map from $M$ to the stack integrating a complex Lie algebroid has nondiscrete fibres.

Consider  $\complex^2$ with coordinates
$z=x+iy$ and $w=s+it$ and the involutive structure  spanned
by $\del / \del x + i\ \del / \del y$ and 
$\del/\del t-(x+iy) \ \del / \del s,$ or, in complex notation, $\del / \del \overline{z}$ and $\del
/ \del t - z \ \del / \del s$.  When $y\neq 0$, this is a complex
structure, while when $y=0$, it contains the real subspace spanned by
$ \del/ \del t - x\ \del / \del s$.  The integrals for this structure
are generated by $z=x+iy$ and $\zeta=s+zt$.  On the complexification
$\complex^2_\complex = \complex ^4$,  $x$, $y$, $s$, and $t$ may 
have complex values, and then the map $(z,\zeta):\complex^2_\complex
\to\complex ^2$ is a submersion whose fibres are the leaves of the
extended foliation; thus, the leaf space (and hence
the stack which integrates the structure) may be identified
with the complex $(z,\zeta)$ plane.

What is singular here is the map $\phi$ from the original $\complex^2 = \reals^4$ to this stack.   When the variables $(x,y,s,t)$ are real, $\phi$ is a local diffeomorphism, except on the hypersurface $y=0$, where each of the orbits of the vector field  
$ \del/ \del t -  x\ \del / \del s$  is mapped to a constant.  The image of this hypersurface is the subset of the $(z,\zeta)$ plane on which the variables are both real, and, as is clearly described by Eastwood and Graham \cite{ea-gr:involutive}, the map $\phi$ realizes the (real) blow-up of $\reals^2$ in $\complex^2$.  

A similar involutive structure was constructed by Lebrun and Mason \cite{le-ma:zoll} on the projectivized complexified tangent bundle of a surface with affine connection; the singular curves in their example are the geodesics.  

\section{Boundary Lie algebroids}
\label{sec-boundary}
This section exhibits CLAs which are
neither involutive systems nor the complexification of real Lie
algebroids.  The example is taken from work in progress by Leichtnam, Tang, and the author \cite{le-ta-we:poisson} on K\"ahler geometry and deformation quantization in the setting of CLAs.  The description of  the integration of these CLAs is not complete.  

Let $X$ be a complex manifold of (complex) dimension $n+1$
 with boundary $M$, and
let $\cale _{M,X}$ be the space of complex vector
fields on $X$  whose values along $M$
lie in the induced CR structure $G_{M,X}$. $\cale _{M,X}$ is a module
over $C^{\infty}(X)$  and is closed under bracket.
The following lemma shows that that it may be identified
with the space of sections of a complex Lie algebroid $E_{M,X}$.

\begin{lemma}
\label{lemma-module}
$\cale _{M,X}$ is a locally free $C^\infty(X)$-module.
\end{lemma}
\pf
Away from the boundary, $\cale _{M,X}$ is the space of sections of $T_\complex M$,
hence locally free.
Near a boundary point, choose a defining function $\psi$, i.e. a function which vanishes on the boundary and has no critical points there.  Next, 
choose a local basis $\vbar_1,\ldots,\vbar_n$ of
$G_{M,X}$ and  extend  it to a linearly independent set of sections of
$T^{0,1}X$, still denoted by $\vbar_j$,
defined in an open subset of $X$, to be shrunk as necessary.
Let $v_j$ be the complex conjugate
of $\vbar_j$.  These vectors all annihilate $\psi$ on $M$;
there is no obstruction to having them annihilate $\psi$ everywhere.
Next, choose a local section $\vbar_0$ of $T^{0,1}X$
such that $\vbar_0\cdot \psi = 1$,  and let $v_0$ be its conjugate.  This
gives a local basis $(v,\vbar)$ for the complex vector fields.
Such a vector field
belongs to $\cale_{M,X}$ if and only if, when it is expanded with respect to
this basis, the coefficients of $\vbar_0$ and all the $v_j$ vanish
along $M$.  Since this means that all these coefficients are
divisible by $\psi$ with smooth quotient,  setting $u'_0=\psi\vbar_0$, $u'_j=\vbar_j$ for
$j=1,\ldots,n$, and $u_j=\psi v_j$ for $j=0,\ldots,n$ produces  a local basis
$(u,u')$ for
$\cale_{M,X}$.
\qed

To integrate the boundary Lie algebroid $E_{M,X}$, assuming analyticity as usual, one may begin by extending $X$ slightly beyond $M$, so that $M$ becomes an embedded hypersurface.  In the complexification $X_\complex$, $M$ extends to a submanifold 
$M_\complex$ of complex codimension one.  The CR structure on $M$ extends (see Section \ref{subsubsec-cr}) to the tangent bundle $E'$ of the CR foliation on $M_\complex$.  The holomorphic continuation of $E_{M,X}$ is then the holomorphic Lie algebroid whose local sections are the vector fields on $X_\complex$ whose restrictions to $M_\complex$ have their values in $E'$.

What is the groupoid of this Lie algebroid over $X_\complex$?  Over the complement of $M_\complex$, the Lie algebroid is the tangent bundle, so the groupoid could be taken to be  the pair groupoid.   Since $M_\complex$ has complex codimension one, though, its complement generally has a nontrivial fundamental group, and the fundamental groupoid or one of its nontrivial quotients might be appropriate as well.  The choice depends in part on compatibility with the choice made on $M_\complex$ itself.  

Over $M_\complex$, the image of the anchor of the extended Lie algebroid is the tangent bundle $E'$ to the CR foliation, but now, unlike in the pure CR situation, there is nontrivial isotropy.  To describe this isotropy, note that, at each point $x$ of $M_\complex$, there is a flag $E'_x \subset T_x M_\complex \subset T_x X_\complex.$  The isotropy algebra may identified with the endomorphisms of the normal space $T_x X_\complex / E'_x$ which vanish on $T_x M_\complex $.   Given two points $x$ and $y$ in 
$M_\complex$, there are morphisms in the integrating groupoid from $x$ to $y$ 
if and only if $x$ and $y$ lie in the same leaf of the CR foliation.  Each such morphism 
is then a linear map 
 $T_x X_\complex / E'_x   \rightarrow T_y X_\complex / E'_y$ whose restriction
 $T_x M_\complex / E'_x   \rightarrow T_y M_\complex / E'_y$ coincides with the linearized holonomy map along any path in the leaf.  (Assume that the complexification is small enough so that the leaves are simply connected.)    In particular, when $x=y$, the isotropy group consists of the automorphisms of $T_x X_\complex / E'_x$ which fix 
 $T_x M_\complex / E'_x.$  (Compare the author's discussion in Section 6 of \cite{we:groupoids}, where the Lie algebroid and its integrating groupoid are studied for the vector fields tangent to the boundary of a real manifold, as well as the treatment by Mazzeo \cite{ma:elliptic} of vector fields tangent to the fibres of a submersion on the boundary.
Finally, a slightly different, class of vector fields on a manifold with fibred boundary is used by Mazzeo  and Melrose \cite{ma-me:pseudodifferential}.)

When $x$ lies on the real hypersurface
$M$,  the space above admits an explicit description in terms of the CR geometry.
Over $M$, 
$T X_\complex$ restricts to $T_\complex X$, $T M_\complex$
is just $T_\complex M$, and $E'$ is the CR structure 
$G_{M,X}=T_\complex M \cap T^{0,1}_J X.$  Thus, the isotropy of the integrating groupoid consists of the automorphisms of $T_\complex X/T_\complex M \cap T^{0,1}_J X$ which fix its codimension one subspace $T_\complex M/T_\complex M \cap T^{0,1}_J X.$  These automorphisms act on the complexified normal bundle $T_\complex X/T_\complex M$, and those which act trivially on the normal bundle are ``shears'' which may be identified with the additive group of linear maps from that normal bundle to 
$T_\complex M/T_\complex M \cap T^{0,1}_J X.$   The choice of a defining function trivializes the normal bundle, so  the isotropy is an extension of the automorphism (or ``dilation") group of the normal bundle by the abelian group 
$T_\complex M/T_\complex M \cap T^{0,1}_J X.$  

The  preceding description of the integrating groupoid is not complete, since it lacks
an explanation of how the piece over the interior and the piece over the boundary fit together. In particular, if one were to use the fundamental groupoid on the interior, as described above, it may be necessary to use a covering of the automorphisms of the line bundle on the boundary.

\section{Generalized complex structures}
In the rapidly developing subject of {\bf generalized geometry}, originated by
Hitchin \cite{hi:generalized}, the tangent bundle $TM$ of a manifold with its
Lie algebroid structure is replaced by the {\bf generalized tangent bundle}
$\calt M$, which is the direct sum $TM\oplus T^*M$ equipped with the
Courant algebroid structure consisting of the bracket
$$
[\![(\xi_{1},\theta_{1}),(\xi_{2},\theta_{2})]\!] =
\left( [\xi_{1},\xi_{2}] , \mathcal{L}_{\xi_{1}}\theta_{2} -
\mathcal{L}_{\xi_{2}}\theta_{1} - \frac{1}{2}d(i_{\xi_{1}}\theta_{2} -
i_{\xi_{2}}\theta_{1})\right)  ,
$$
the anchor $\calt M \to TM$ which projects to the first summand, and
the symmetric bilinear form

\[
\langle(\xi_1 ,\theta_1 ),(\xi_2 ,\theta_2 )\rangle = 
\frac{1}{2}(i_{\xi_1}\theta_2 + i_{\xi_2}\theta_1 ).
\]

Like the tangent bundle, $\calt M$ may be complexified to the ``complex
Courant algebroid'' $\calt_\complex M$.  It is not a complex Lie algebroid,
but it contains many CLAs, in particular the {\bf complex Dirac
structures}, i.e. the (complex) subbundles $E$ which are maximal isotropic
for the symmetric form and whose sections are closed under the
bracket.  For instance, if $A\subseteq T_\complex M$ is an involutive
system and $A^\perp \subseteq T^*_\complex M$ is its annihilator, then
$A\oplus A^\perp$ is a complex Dirac structure.  

Of special interest among the complex Dirac structures are
those for which $E\oplus \overline{E} = \calt M$.  These are called
{\bf generalized complex structures} and are the $-i$ eigenspaces of
(the complexifications of) integrable almost complex structures
$\calj:\calt M\to\calt M$; the integrability condition here is that the 
Nijenhuis torsion is zero, the usual bracket of vector
fields in the definition of the torsion being replaced by the Courant bracket.

In particular given a complex structure $J:TM\to TM$, with associated
CLA $T^{0,1}_J M$, the direct sum with its annihilator is the
generalized complex structure $\calt^{0,1}_J M=T^{0,1}_J M \oplus
{T^{1,0}_J}^* M$.  The image of the anchor is the involutive system
$T^{0,1}_J M$, but $\calt^{0,1}_J M$ itself is not an involutive
system, since the kernel of its anchor is the nontrivial bundle
${T^{1,0}_J}^* M, $.  Also,  $\calt^{0,1}_J M$ is not isomorphic to
the complexification of a real Lie algebroid, since the image of its
anchor is not invariant under complex conjugation.

Another kind of example arises from
symplectic structures on $M$, viewed as 
bundle maps $\omega:TM\to T^*M$.  Here, 
the generalized complex structure $E_\omega$ is defined to be the graph of the
complex 2-form $i\omega$.  This time, the anchor is bijective, so, as
a Lie algebroid, $E_\omega$ is isomorphic to $T_\complex M$.  

What is the integration, in the sense of this paper, of a generalized complex structure?   First, let $J$ be a complex structure on $M$, 
$\calt^{0,1}_J M=T^{0,1}_J M \oplus {T^{1,0}_J}^* M$ the corresponding generalized structure.  Complexifying $M$ and $J$ as in Section \ref{subsubsec-complex} gives a foliation on $M_\complex$.  The groupoid which integrates the holomorphic continuation of 
$\calt^{0,1}_J M$ is the semidirect product  groupoid obtained from the action of the holonomy groupoid of the foliation (via the ``Bott connection") on its conormal bundle.  (This is just the holomorphic version of a construction  by Bursztyn, Crainic, Zhu, and the author \cite{bu-cr-we-zh:integration}.) This action groupoid is equivalent to the holomorphic leaf space $M$  carrying the
cotangent bundle  ${T^{1,0}_J}^* M$ of additive groups as its isotropy.  The corresponding stack is a the bundle over $M$ whose fibres are the ``universal classifying stacks" of the cotangent spaces.  

Next let $\omega$ be a symplectic structure on $M$.   Since the generalized complex structure $E_\omega$ is isomorphic to $T_\complex M$, its integration must be that of $T_\complex M$, i.e. the holomorphic point, perhaps carrying the fundamental group of $M$ as isotropy.    To see what has become of $\omega$, it is best to look again at (real and complex) Dirac structures.  

As a subbundle of $\calt M$, a Dirac structure $E$ carries a natural skewsymmetric bilinear form, the restriction of 
\[
B(\xi_1 ,\theta_1 ),(\xi_2 ,\theta_2 )) = 
(1/2)(i_{\xi_1}\theta_2 - i_{\xi_2}\theta_1 ).
\] 
It is shown in \cite{bu-cr-we-zh:integration} that this form gives rise to a multiplicative closed 2-form on a groupoid integrating $E$, producing a {\bf presymplectic groupoid}.     
Applying this construction to the holomorphic extension of any complex Dirac structure
$E$ shows that its integration as a CLA  is a holomorphic symplectic groupoid over $M_\complex$.  In particular, for $E_\omega$ or any other complex Poisson structure, it is a holomorphic symplectic groupoid.  For $E_J$, or any other direct sum of an involutive structure with its annihilator, the restriction of $B$ is zero, and hence so is the presymplectic structure on the integrating groupoid.

\section{Further topics and questions}
A notion of integration for complex Lie algebroids has been proposed in this paper.  There are many interesting questions  about other extensions of Lie algebroid theory to the complex case, including the relation between these extensions and the integration construction proposed here.   Some examples conclude this paper.

\subsection{Integrability}
Does the integrability criterion of Crainic and Fernandes \cite{cr-fe:lie} apply in the holomorphic case?  What are the conditions on an analytic CLA which determine whether its holomorphic continuation is integrable?
What can one do in the nonanalytic case?

\subsection{Cohomology}
 A ``van Est" theorem of Crainic \cite{cr:differentiable} describes the relation between the cohomology of a Lie algebroid and that of its integrating groupoids.  The definition of cohomology extends in a straightforward to CLAs (for instance, it gives the Dolbeault cohomology in the case of a complex structure).  Is there a  van Est theorem in this case, too?

\subsection{Bisections}
One consequence of the integration of a Lie algebroid $E$ is that the submanifolds of an integrating groupoid which are sections for the source and target maps form a group whose Lie algebra in some formal sense is the space of sections of $E$.  Is there a similar construction for the case of a complex Lie algebroid?  Some hints might come from  the constructions by Neretin \cite{ne:complex} and Segal \cite{se:definition} (also see Yuriev \cite{yu:infinite}) of a semigroup which in some sense integrates the complexified Lie algebra of vector fields on a circle.  Conversely, a general construction for CLAs could provide complexifications for the diffeomorphism groups of other manifolds. 
 
\subsection{Quantization}
Once a Lie algebroid $E$ has been integrated, the groupoid algebra of an integrating groupoid may be considered, following Landsman and Ramazan \cite{la-ra:quantization}, as a deformation quantization of the Poisson structure on the dual bundle $E^*$, or as a completion of Rinehart's \cite{ri:differential} universal enveloping algebra  of $E$. Is there a  corresponding application for the integration of a CLA?

On the other hand, given a complex Poisson structure $\Pi$ on $M$, it defines a CLA structure on the complexified cotangent bundle.  Integration of this structure should give a holomorphic symplectic groupoid which should be somehow related to the deformation quantization of $(M,\pi)$.  On the formal level (without integration), it is possible \cite{le-ta-we:poisson} to extend the methods of Karabegov \cite{ka:separation} and 
Nest and Tsygan \cite{ne-ts:deformations} to construct deformation quantizations of certain boundary Lie algebroids as in Section \ref{sec-boundary} above.

\subsection{Connections and representations}
If $E$ is a CLA over $M$ and $V$ is a complex vector bundle $V$,  an $E$-connection on $V$ is a map $a\mapsto \nabla_a$ from the  sections of $E$ to the $\complex$-endomorphisms of the sections of $V$ which satisfies the conditions $\nabla_{fa} u = f\nabla_a u$ and $\nabla_a gu = g \nabla_a u + (\rho(a)g)u$.    The connection is flat and is also called a representation of $E$ on $V$ if the map $a\mapsto \nabla_a$ is a Lie algebra homomorphism.

For instance, if $E$ is  a complex structure, a representation of $E$ on $V$ is a holomorphic structure on $V$.   More generally, representations of CR structures correspond to CR vector bundles, as in the work of Webster \cite{we:integrability}.  After complex extension, an analytic  representation of an analytic CR structure becomes a flat connection along the leaves of the CR foliation, which leads to a holomorphic vector bundle on the complexification.  

If $E$ is the generalized complex structure associated to a complex structure, a representation is a holomorphic structure together with a holomorphic action of the holomorphic cotangent bundle by bundle endomorphisms of the representation space.

\subsection{The modular class}
The modular class of a Lie algebroid, introducted by Evans, Lu, and the author \cite{ev-lu-we:transverse} is the obstruction to the existence of an ``invariant measure."   Its definition extends directly to the case of CLAs.  For a complex structure, the modular class is the obstruction to the existence of a Calabi-Yau structure.

\end{document}